\documentclass[11pt]{article}
\usepackage{amsmath}
\usepackage{amssymb}
\usepackage{url}
\newcommand\sLP{\\[\smallskipamount]}
\newcommand\sPP{\\[\smallskipamount]\indent}
\newcommand\mLP{\\[\medskipamount]}
\newcommand\mPP{\\[\medskipamount]\indent}
\newcommand\bLP{\\[\bigskipamount]}
\newcommand\bPP{\\[\bigskipamount]\indent}
\newcommand\RR{\mathbb{R}}
\newcommand\Ga\Gamma
\newcommand\La\Lambda
\newcommand\al\alpha
\newcommand\be\beta
\newcommand\ka\kappa
\newcommand\la\lambda

\newcommand\iy\infty
\numberwithin{equation}{section}
\begin{document}
\title{Josef Meixner: his life and his orthogonal polynomials}
\author{Paul L. Butzer and Tom H. Koornwinder}
\date{}
\maketitle
\begin{abstract}
This paper starts with a biographical sketch of the life of Josef
Meixner. Then his motivations to work on orthogonal polynomials
and special functions are reviewed. Meixner's 1934 paper introducing
the Meixner and Meixner-Pollaczek polynomials is discussed in detail.
Truksa's forgotten 1931 paper, which already contains the Meixner
polynomials, is mentioned. The paper ends with a survey of
the reception of Meixner's 1934 paper.
\end{abstract}
\section{Introduction}
This paper grew from the idea that certain names occurring in
the {\em Askey scheme of hypergeometric orthogonal polynomials}
\cite[Appendix]{33}, \cite[p.184]{1}
deserve more historical explanation, both biographically and
mathematically. We decided to start with Meixner, whose name occurs
in the Meixner polynomials and the Meixner-Pollaczek polynomials.
The German theoretical physicist Josef Meixner lived from 1908 until
1994. He spent most of his career in Aachen.

Meixner's paper \cite{Me2}
containing the polynomials named after him appeared in 1934 as one
of his first papers, while he was still in M\"unchen.
In current terminology this paper introduced Sheffer polynomials
by giving two equivalent definitions, one in terms of generating
functions, and then classified all orthogonal polynomials which are
Sheffer polynomials. Two new classes came up in this classification:
the Meixner and Meixner-Pollaczek polynomials. Reception of this
work was not quick, but after some decades the paper got the
recognition it deserved. 

Meixner did much further work on special functions, in particular
on spheroidal functions and Mathieu functions, on which he published
with Sch\"afke the research monograph \cite{B}. In physics Meixner
is best known for his work on the thermodynamics of irreversible
processes, see his article \cite{MeixRe} with Reik in
{\em Handbuch der Physik}. In this area he was competing, in a certain
sense, with Nobel prize winner Prigogine.

The contents of this paper are as follows.
In Section \ref{33} we describe Meixner's life.
Section \ref{1} contains quotations about his motivations to
work on orthogonal polynomials and special functions.
Next, in Section \ref{21}, we give a detailed sketch of his
paper \cite{Me2} introducing the Meixner polynomials.
However, there was a predecessor on Meixner polynomials:
Truksa \cite{15}. This is discussed in Section \ref{20},
where we also pay attention to Pollaczek's 1950 paper \cite{25},
which rediscovered the Meixner-Pollaczek polynomials.
Finally, in Section \ref{34}, we review the reception of Meixner's
1934 paper \cite{Me2} over a period of more than 80 years.
A concluding remark mentions the irony of history in fixing names
for classes of polynomials and for theorems.

We start the References with a selection of Meixner's publications:
\cite{Me1}--\cite{Me16} give papers on orthogonal polynomials and
special functions,
\cite{Me7}--\cite{Me15b} papers on
spheroidal wave functions and Mathieu functions,
and \cite{B}--\cite{MeSchWo} books and Handbuch articles.
\section{Biographical notes}
\label{33}
Josef Meixner, born 24 April 1908 in Percha (now part of Starnberg,
Bayern), studied mathematics, physics and physical
chemistry at the {\em Universit\"at M\"unchen} during 1926--1931.
Shortly after receiving
the upper level teaching certificate in February 1931 he obtained
his doctorate under Arnold Sommerfeld\footnote{Arnold Sommerfeld (1868--1951), the internationally respected
theoretical physicist, was counted
``among the few German scientists who were untainted 
with regard to Nazi affiliation'', see Eckert \cite{H} and
\url{http://www.encyclopedia.com/topic/Arnold_Johannes_Wilhelm_Sommerfeld.aspx}\,.} with the thesis
{\em Die Greensche Funktion des wellenmechanischen \mbox{Keplerproblems}}.
It was
the first treatment of Green's function in quantum mechanics and
the first example of the $S$-Matrix Theory of Wheeler (1937) and  Heisenberg (1943).
Ten pages of this thesis were devoted to mainly new
investigations of confluent hypergeometric series needed, series
which occur rather frequently in a good part of Meixner's
oeuvre.

He was an assistant at the Institute of Theoretical Physics in M\"unchen
until 1934 when he became an assistant to
Karl Bechert\footnote{Karl Bechert (1901--1981), a Sommerfeld student and Professor of Theoretical Physics
at the {\em Universit\"at Gie\ss en} since 1933, was appointed Rector
of this university by the Americans
in 1945/46,  he being one of the few
who was not politically polluted.
Thereafter he became Director of the
Institute of Theoretical Physics at the  {\em Universit\"at Mainz}, refounded
by the French military government in 1946.
See \url{http://www.regionalgeschichte.net/bibliothek/texte/biographien/bechert-karl.html}\,.}
at the
{\em Universit\"at Gie\ss en}. In 1934 he also became a member of the SA and in 1937
of the NSDAP (see \cite[p.325, footnote 1]{C}).
In Gie\ss en he received the
\emph{Habilitation} degree together with a teaching assignment in
theoretical physics in 1936/37. After lecturing at the
universities of Marburg and Berlin in the period
1938--1942\footnote{These were the years that Meixner was also
quite active
as reviewer for {\em JFM} (do structured search in the reviewer field in
\url{https://zbmath.org/}) on physics related papers; earlier he had reviewed a few papers on orthogonal polynomials for
{\em Zentralblatt}.}, he was appointed
as an extraordinary Professor of
Theoretical Physics at the
{\em Technische Hochschule Aachen} in
December 1942, as successor to Prof.\ Wilhelm Seitz (1872--1945).
But he could not take up this position as he was not released from
the Armed Forces; he was stationed at the weather station at
Vads\"o, Kirkenes, Norway, since September 1941. Released from the
forces in summer 1943, he was sent by Major Gottfried Eckart
(Heidelberg), in charge of high frequency research, to
Murnau (south of M\"unchen) to do this research together with F. Sauter and
A.~W. Maue, who were also former students of Sommerfeld\footnote{This move
was possibly also intended as giving a
counterweight to Sommerfeld's successor in office, Wilhelm
M\"uller (born 1880), known to be a notorious Nazi.}. It was here
that Meixner began his work on Mathieu and spheroidal
functions.

Meixner was welcomed by Aachen's (provisional) Rector, Prof.\ Paul
R\"ontgen\footnote{As to R\"ontgen: ``during the Nazi era he
remained neutral as far as possible $\ldots$ non-Nazi elements speak
of him as being one of the few respectable and acceptable members
of the TH staff $\ldots$ during the Nazi regime as well as now
$\ldots$'' \cite{C}. See also \cite{e} about the TH Aachen during the
war.} (1881--1968, born in Aachen) on 30 September 1945.
However, his wartime appointment in Aachen could not be automatically
continued, but for reemployment a {\em Persilschein} (denazification
certificate) was needed which was written by
Sommerfeld\footnote{14 December 1945,
Deutsches Museum, Archiv, NL 89, 020, Mappe 8,3;
included in questionnaire of the British Army of 30 July 1946, see
\cite[p.325, footnote 1]{C}.}. He writes that
he has known Meixner for 17 years, that Meixner has never been a supporter
of the Nazi system, but that 
in Meixner's given circumstances in Gie\ss en in 1934
it would have been very difficult for him
to avoid membership of the SA.
He urges  TH Aachen to try its best
to keep ``such an outstanding teacher and scholar'' as Meixner in office.
Meixner was appointed as an ordinary Professor on a 	
personal title in 1948, and as a
regular chairholder in 1951. In fact, he was not only a member of the
Faculty of Natural Sciences, but later also of the Faculty of
Electrical Engineering. The only physicist at the time in Aachen was
Prof.\ Wilhelm Fucks (1902--1990), one of the politically active
National Socialists at the TH (\cite{C} pp.~320--323). Meixner stayed
in Aachen despite several attractive offers from other
universities in Germany and abroad. He was visiting or research
professor in the USA at New York University (1954), Michigan State University
(1955/56), University of Michigan (1961), University of Southern
California (1964), Brown University (1965), Lehigh University (1969/70)
and in Japan at Hokkaido University (1975).
He received an
honorary doctorate from the {\em Universit\"at K\"oln} in 1968.
He was a Fellow of the
{\em Nordrhein-Westf\"alische Akademie der Wissenschaften und der K\"unste}.
He is the author of 150 publications.

In 1954 Meixner wrote together with the mathematician Friedrich
Sch\"afke the well-known book \cite{MeSch} on Mathieu functions
and spheroidal functions. In 1956 he published an article on special
functions in the {\em Handbuch der Physik} \cite{Meix}.

A major part of Meixner's research work concerned the
thermodynamics of irreversible processes and he is counted as one
of the founding fathers of that field. The fundamental principles
of this theory he published together with G.~Reik, in their famous
article in the {\em Handbuch der Physik} \cite{MeixRe} in 1959.
The Nobel Prize for Chemistry for the year 1977
\footnote{See
\url{http://www.nobelprize.org/nobel_prizes/chemistry/laureates/1977/}\,.} was awarded to
Ilya Prigogine (1917--2003), see \cite{G}, especially for his work on
the thermodynamics of irreversible processes. It was said the
Prize Committee spent 4 overtime hours before reaching its
decision, and there were rumours that Meixner was also a candidate.
Meixner's first basic paper in the matter dates back to 1941,
Prigogine's work started in 1947. In fact, D. Bedeaux and I. Oppenheim
\cite{h} in their obituary notice on Mazur write: ``Josef Meixner
in 1941 and, independently, Ilya Prigogine in 1947 set up a
consistent phenomenological theory of irreversible processes,
incorporating both Onsager's\footnote{Lars Onsager (1903--1976)
received the Nobel Prize for Chemistry in 1968. For a paper primarily
based on an interview with Onsager shortly before his death see
\cite{F}.} reciprocity theorem and the explicit calculation for
some systems of the so-called entropy source strength. Shortly
thereafter, Mazur and de Groot joined this group as founding
fathers of the new field of nonequilibrium thermodynamics''.
Furthermore, A.~R. Vasconcellos et al.\ \cite{E} write:
``J. Meixner, over twenty years ago
in papers that did not obtain a deserved diffusion  gave some
convincing arguments to show that it is unlikely that a
nonequilibrium  $\ldots$''. In 1975 both Prigogine and Meixner lectured
at an Academy session \cite{PrigMeix} in D\"usseldorf.
See also \cite{D}. Another appraisal of Meixner's work on
thermodynamics appeared in \cite{f}.

After Sommerfeld's death in 1951 Meixner and F. Bopp were
entrusted by his last will to complete his manuscript on
thermodynamics and statistics \cite{B}. ``Especially Meixner 's
section on thermodynamics of irreversible
process fits in admirably in Sommerfeld's masterpiece.''
(see~\cite{b} and \cite{a}).

We could trace three doctorates with Meixner as {\em Referent}
(first advisor): Ingo M\"uller\footnote{See
the Mathematics Genealogy Project,
\url{http://www.genealogy.math.ndsu.nodak.edu/id.php?id=48491};
Meixner was only a {\em Korreferent} for the other three mentioned there.}
(born 1936, doctorate in 1966), Wolfgang Kern and
Gunter Weiner\footnote{W. Kern confirmed to us that Weiner and he
were PhD students of Butzer.}.
Kern and Weiner, together with Meixner, are the authors of the booklet \cite{KeWeMeix}.

Meixner retired from his professorship in Aachen in 1974.
He passed away in 1994.
From an obituary note by Schl\"ogl \cite{c}, \cite{d} we quote
the following.
``Josef Meixner applied severe but just standards in both research
and teaching. His well-founded, not-always-comfortable judgments were
highly respected by his colleagues.''
\section{Meixner's work in orthogonal polynomials and special functions:
his motivations}
\label{1}
One of us (PB) had access to Meixner's handwritten {\em Erinnerungen},
which were later typewritten (44 pp.).
We quote a few parts of this document which give insight about
Meixner's motivations for his work on orthogonal polynomials and
special functions.
\mLP
``The Hermite polynomials\footnote{Meixner meant here the Hermite
polynomials orthogonal with respect to the weight function $\exp(-x^2/2)$;
these are nowadays often notated by $He_n(x)$.}
have a generating function of the form
\begin{equation*}
f(t) \exp(xt) = \sum_{n=0}^{\infty} H_n (x)\,\frac{t^n}{n!}\,.
\end{equation*}
The mathematician Salomon Bochner\footnote{Meixner did not write
a joint paper with Bochner but in a footnote to his paper in
J.~London Math.\ Soc.\ \cite{Me2} he writes that in the simpler
particular case $u(t)=t$ Bochner was involved in its investigation
with him.} and I pursued the question as to what functions $f(t)$
lead to other orthogonal polynomial systems. The result was: only
the Hermite polynomials do have such a generating function.
I broadened the question and sought
for polynomial systems having a generating function
\begin{equation}
f(t) \exp{(x\,u(t))} = \sum_{n=0}^{\infty} p_n(x)\,\frac{t^n}{n!}\,,
\label{2}
\end{equation}
where $f(t)$ and $u(t)$ are formal power series in $t$, with
$f(0)=1$, $u(0)=0$ and $u'(0)=1$. The solution of this
problem leads to interesting results (published 1934). Apart from
the Hermite, Laguerre and Charlier polynomials there turned out to be
two further classes of polynomials, one having discrete weights
and one with a continuous weight distribution.''
\mLP
Later on he writes:
\mLP
``These polynomial classes were rediscovered by Pollaczek and
it lasted some years until someone realized that these
polynomials were found by me. In the meantime they are called
Meixner's polynomials, and especially the polynomial class having
discrete weights has found numerous applications, see e.g. Chihara
\cite{A} $\ldots$''.
\\[\medskipamount]
``With Bochner [1899--1982], who came circa 1929 to Munich, I was a close
friend and I owe him not only much from the mathematical
side; also his literary interests did indeed stimulate me.
It was perhaps his most fruitful time, at the beginning of true
functional analysis. His book Fourier Integrals, which in a hidden form
already contains an essential part of the theory of distributions,
and which also contains many other of his own results, testifies
this.''{\footnote{Conversely, Bochner valued Meixner's work. Indeed,
in October
1977 Meixner was invited to the birthday meeting of Prof.\ J.~S. Frame
who there recalled Bochner's lecture at Cambridge, UK, in 1933, in which
he had placed special emphasis on Meixner's work on orthogonal polynomials
having a generating function of special type; i.e. Meixner's paper 
\cite{Me2}.}}
\mLP
Elsewhere Meixner writes:
\mLP
``Bochner foresaw the coming political development
very clearly, and I recall when we, surely at the end of 1932, stood before a
bulletin board of the {\em V\"olkischer Beobachter}
and he said: `Now it is 
almost time
that I must depart'. When I [at age 24] replied that then I would also like  
to leave,
he replied: You remain here; nothing will happen to you and for us there
are too few places in the world.''
\mLP
On p.40 of his {\em Erinnerungen} Meixner writes:
\mLP
``The special functions of
mathematical physics especially interested me already since my 
dissertation. In the
last war year I occupied myself with a possible application of
Mathieu functions and the spheroidal functions which are
considerably more difficult to handle than most of the special
functions studied until then. There existed many individual
results but not even a logically consistent nomenclature. My
contribution involved a rational terminology and the derivation of
numerous new properties over a period of more than ten years. In
this matter my later collaborator
F.~W. Sch\"afke\footnote{F.~W. Sch\"afke in his letter of April 1978
addressed to Meixner writes that since at the end of the war he
had neither a scientific advisor nor belonged to a ``school'' (his
advisor Harald Geppert had died 1945) his formative years were
essentially influenced through the collaboration with Meixner in
their work on Mathieu functions. It formed the beginning of his
own work of special functions $\ldots$ . Finally Sch\"afke adds: ``You
were for me fatherly-kind and encouraging, a shining scientific
and human example.''} supported me and our association
culminated in 1954 in a book (see \cite{MeSch}) concerning these
functions which represent their
first systematic theory.''
\mPP
Meixner made great impression by his lecture at the International
Christoffel Symposium (Aachen, Monschau, November 1979), which
was organized by one of us (PB) together with F.~Feh\'er.
The resulting paper \cite{Me16} (1981), one of his last publications, discusses
the {\em Christoffel-Darboux formula} \cite[Theorem 3.2.2]{8}:
the way it was first derived by Christoffel (who was prior to Darboux)
and its usage, in particular for approximation of functions by their
partial sum expansions in terms of orthogonal polynomials.
This paper, with its precise and stimulating presentation, 
is a nice testimony of Meixner's
lasting interest in orthogonal polynomials
and his keeping up with recent results in literature.
\section{Meixner's introduction of the Meixner and
Meixner-Pollaczek polynomials}
\label{21}
As could already be read in Section \ref{1}, Meixner introduced the
orthogonal polynomials named after him in his 1934 paper \cite{Me2}
and this work arose from his contact with Bochner on a special case
of the problem considered in \cite{Me2}. Actually, this purely
mathematical paper \cite{Me2}
was Meixner's first published paper after his
paper \cite{Me1}, which was based on his more physically oriented
dissertation.

Meixner classifies in \cite{Me2} all orthogonal polynomial systems
$\{p_n\}_{n=0,1,2,\ldots}$ which have a generating function
\eqref{2} with $f$ and $u$
as specified after \eqref{2}. The results of the classification are
the cases I--V in \cite[p.~11]{Me2}.
We want to observe that, on the one hand, \eqref{2} can be written
more generally for not necessarily monic polynomials $p_n$ as
\begin{equation*}
f(t) \exp{(x\,u(t))} = \sum_{n=0}^{\infty} c_n\,p_n(x)\,t^n,
\end{equation*}
for certain real non-zero coefficients $c_n$.
With dilation of the argument $x$ we can also replace the condition
$u'(0)=1$ by $u'(0)\ne0$. This larger flexibility is useful when
considering examples.
On the other hand, with suitable translation
of the argument $x$, we can make for \eqref{2} the further requirement
that $f'(t)=0$ or, equivalently, that $p_1(x)=x$. This is used
by Meixner in the proof of his classification.

With names and notation as
currently used (see \cite[Chapter 9]{1}), these five cases read as follows.
\paragraph{Case I}
Hermite polynomials $H_n(x)$, \cite[Section 9.15]{1}:
\begin{equation}
\exp(xt)\,e^{-\frac14 t^2}
=\sum_{n=0}^\iy\frac{H_n(x)}{2^n\,n!}\,t^n,
\label{3}
\end{equation}
orthogonal on $(-\iy,\iy)$ with respect to the weight function
$e^{-x^2}$.
\paragraph{Case II}
Laguerre polynomials $L_n^{(\al)}(x)$, \cite[Section 9.12]{1}:
\begin{equation*}
\exp\left(\frac{xt}{t+1}\right)(1+t)^{-\al-1}
=\sum_{n=0}^\iy(-1)^n\,L_n^{(\al)}(x)\,t^n,
\end{equation*}
orthogonal on $[0,\iy)$ with respect to the weight function
$x^\al e^{-x}$ ($\al>-1$).
\paragraph{Case III}
Charlier polynomials $C_n(x;a)$, \cite[Section 9.14]{1}:
\begin{equation*}
\exp\big(x\log(1+t)\big)\,e^{-at}
=\sum_{n=0}^\iy\frac{(-a)^n\,C_n(x;a)}{n!}\,t^n,
\end{equation*}
orthogonal on $\{0,1,2,\ldots\}$ with respect to the weights
$a^x/x!$ ($a>0$).
\paragraph{Case IV}
Meixner polynomials $M_n(x;\be,c)$, \cite[Section 9.10]{1}:
\begin{equation*}
\exp\left(x\log\left(\frac{1-c+t}{1-c+ct}\right)\right)
\left(1+\frac{ct}{1-c}\right)^{-\be}
=\sum_{n=0}^\iy\frac{(\be)_n}{n!}\left(\frac c{c-1}\right)^n
M_n(x;\be,c)\,t^n,
\end{equation*}
orthogonal on $\{0,1,2,\ldots\}$ with respect to the weights
$c^x(\be)_x/x!$ ($\be>0$, $0<c<1$).
\paragraph{Case V}
Meixner-Pollaczek polynomials $P_n^{(\la)}(x;\phi)$, \cite[Section 9.7]{1}:
\begin{equation*}
\exp\left(ix\log\left(\frac{2\sin\phi-e^{i\phi}t}
{2\sin\phi-e^{-i\phi}t}\right)\right)
\left(\frac{4\sin^2\phi-4t\sin\phi\cos\phi+t^2}
{4\sin^2\phi}\right)^{-\la}
=\sum_{n=0}^\iy\frac{P_n^{(\la)}(x;\phi)}{(2\sin\phi)^n}\,t^n,
\end{equation*}
orthogonal on $(-\iy,\iy)$ with respect to the weight function
$e^{(2\phi-\pi)x} |\Ga(\la+ix)|^2$\\
($\la>0$, $0<\phi<\pi$).
\bPP
The cases I and II were well-known and Meixner refers for III to
Charlier \cite{4} (1905), Jordan \cite{5} (1926) and
Pollaczek-Geiringer \cite{6} (1928).
But Meixner did not give citations for IV and~V. Still class IV had
already appeared in a paper in 1931 by Truksa \cite{15}, see
Section \ref{20}.

Meixner's starting point in \cite{Me2} is a linear operator $\La$
on the space of polynomials such that
\begin{equation}
D\La=\La t(D),
\label{5}
\end{equation}
where $D$ is the differentiation operator and $t(D)$ is a formal
power series in $D$ of the form
\begin{equation}
t(D)=D+a_2D^2+a_3D^3+\cdots\;.
\label{18}
\end{equation}
Meixner mentions as examples, for the case $t(D)=D$, the operators
\begin{equation*}
\La P(x):=\sum_{\nu=1}^p c_\nu P(x+d_\nu),
\end{equation*}
considered in 1927 by Bochner \cite{2}, and their continuous analogues
\begin{equation}
\La P(x):=\int_{-\iy}^\iy P(x+y)\,d\psi(y).
\label{4}
\end{equation}
Note the case $d\psi(y)=(2\pi)^{-1/2}e^{-y^2/2}\,dy$ of \eqref{4}. Then
\begin{equation*}
\La P(x)=(2\pi)^{-1/2}\int_{-\iy}^\iy P(x+y)\,e^{-y^2/2}\,dy.
\end{equation*}
With this choice of $\La$ we see,
by the generating function \eqref{3} for Hermite polynomials, that
\begin{equation*}
\La P_n(x)=x^n\quad\mbox{if}\quad P_n(x)=H_n(2^{-1/2}x).
\end{equation*}

Meixner observes about the operator $\La$ in the general case that there
is a nonnegative integer $\mu$ such that $\La$ sends any polynomial
of degree $n$ to a polynomial of degree $n-\mu$ if $n\ge\mu$ and to 0
if $n<\mu$. Normalize $\La$ such that $\La x^\mu=1$.
Let $u(s)$ be the formal power series
\begin{equation}
u(s)=s+b_2s^2+b_3s^3+\cdots\quad\mbox{such that}\quad
t(u(s))=s.
\label{19}
\end{equation}
Then, by \eqref{5},
\begin{equation*}
\frac d{dx}(\La e^{xu(s)})=s\La e^{xu(s)},
\end{equation*}
which can be solved as
\begin{equation}
\La e^{xu(s)}=\frac{s^\mu e^{xs}}{\mu!\,f(s)}
\label{6}
\end{equation}
with $f(s)$ a formal power series with constant term equal to 1.
We can expand
\begin{equation}
f(s) e^{xu(s)}=\sum_{n=0}^\iy \frac{P_n(x)}{n!}\,s^n
\label{7}
\end{equation}
with $P_n(x)$ a monic polynomial of degree $n$.
Then, by \eqref{6} and \eqref{7},
\begin{equation}
\La P_n(x)=
\begin{cases}
\binom n\mu x^{n-\mu}&\mbox{if $n\ge\mu$,}\sLP
\;0&\mbox{otherwise.}
\end{cases}
\label{10}
\end{equation}

We have seen that the operator $\La$, implying $t(D)$ and $\mu$,
determines $u(s)$ and $f(s)$.
Conversely, $u(s)$ and $f(s)$ determine the
$P_n(x)$ and next, together with $\mu$, also $\La$.

From now on assume that $\mu=0$, so \eqref{10} takes the simple form
\begin{equation*}
\La P_n(x)=x^n.
\end{equation*}
In combination with \eqref{5} this shows that
\begin{equation}
t(D) P_n(x)=nP_{n-1}(x).
\label{11}
\end{equation}
Also observe from \eqref{7} that
\begin{equation}
f(s)=\sum_{n=0}^\iy \frac{P_n(0)}{n!}\,s^n.
\label{12}
\end{equation}
Meixner states that for a sequence of monic polynomials
$P_n(x)$ the properties \eqref{7} and \eqref{11} are equivalent,
where $t(D)$ of the form \eqref{18} is related to $u(s)$ by \eqref{19}
and $f(s)$ is given by \eqref{12}. Such polynomials $P_n(x)$ are
nowadays called {\em Sheffer polynomials} because of Sheffer's paper
\cite{7} (1939), see the next Section.

Next, in order to arrive at the classification I--V, Meixner assumes
that the $P_n$ are monic orthogonal polynomials,
\begin{equation}
\int_{-\iy}^\iy P_m(x)\,P_n(x)\,d\psi(x)=0\qquad(m\ne n).
\label{9}
\end{equation}
Thus they have to
satisfy a three-term recurrence relation
\begin{equation}
P_{n+1}(x)=(x+l_{n+1})P_n(x)+k_{n+1}P_{n-1}(x)
\label{8}
\end{equation}
with $P_{-1}(x):=0$ and with $k_{n+1}<0$ for $n\ge1$.
Meixner observes that conversely \eqref{8} together
with the negativity of the $k_{n+1}$ implies the orthogonality
\eqref{9} for some monotonic function $\psi$ with infinitely many
jumps and all moments existing. For this result he refers
to Perron's book \cite[p.~376]{3} (1929).
However, note that Meixner's statement cannot be found there explicitly
One has to combine several ingredients in Perron's book
in order to get the result.
In this way Meixner anticipated the formulation
of what is usually called Favard's theorem, as given in
Favard's 1935 paper \cite{36}.
Although Favard was criticized that his paper
gave a well-known result (see for instance W. Hahn's
review of \cite{36} in JFM 61.0288.01), he was also praised that
he was the first to formulate the results in terms of orthogonal
polynomials, and thus Meixner can already be praised for this.

From \eqref{11} applied to \eqref{8} Meixner derives
that the coefficients in \eqref{8} have the form
\begin{equation}
l_{n+1}=l_1+n\la,\qquad
k_{n+1}=n(k_2+(n-1)\ka)\quad(\la\in\RR,\;k_2<0,\;\ka\le0),
\label{14}
\end{equation}
and that $t(u)$ satisfies the differential equation
\begin{equation}
t'(u)=1-\la t(u)-\ka t(u)^2.
\label{13}
\end{equation}
In \eqref{8} with coefficients given by \eqref{14} $l_1$ can be put
equal to zero without essential loss of generality
(by suitable translation of the argument $x$). So we
can work with the three-term recurrence relation
\begin{equation}
P_{n+1}(x)=(x+n\la)P_n(x)+n(k_2+(n-1)\ka)P_{n-1}(x).
\label{15}
\end{equation}
From \eqref{12} together with \eqref{15} Meixner derives the
differential equation
\begin{equation}
\frac{f'(t)}{f(t)}=\frac{k_2t}{1-\la t-\ka t^2}\,.
\label{16}
\end{equation}

Next he factorizes $1-\la t-\ka t^2=(1-\al t)(1-\be t)$, he distinguishes
five different possibilities for $\al,\be$ and he solves \eqref{13}
and \eqref{16} for each of these possibilities.
Thus he arrives at the five cases of the classification, for which he now
explicitly has the functions $t(u)$, $u(t)$ (by functional inversion) and
$f(t)$, and hence the generating function \eqref{2}.
In the case of generic $\al,\be$ he obtains
\begin{equation}
t(D)=\frac{e^{(\al-\be)D}-1}{e^{(\al-\be)D}-\be}\,.
\label{17}
\end{equation}
This degenerates for $\be\to\al$ to $t(D)=\frac D{1+\al D}$, while
in the generic case we observe that
\begin{equation*}
\big(e^{(\al-\be)D}-1\big)g(x)=g(x+\al-\be)-g(x).
\end{equation*}

He obtains the orthogonality measure $d\psi(x)$ by observing
from \eqref{2} and \eqref{9} that
\begin{equation*}
\int_{-\iy}^\iy e^{xu}\,d\psi(x)=
\frac1{f(t(u))}\int_{-\iy}^\iy d\psi(x),
\end{equation*}
and next doing Fourier or Laplace inversion.

Meixner is also able to derive the second order differential or
difference eigenvalue equation satisfied by the polynomials $P_n(x)$.
For this purpose he applies $t(D)$ twice to \eqref{15} in order
to obtain by \eqref{11} that
\begin{equation*}
(n+2)P_n(x)=\big(x+(n+1)(\al+\be)\big)t(D) P_n(x)
+2t'(D) P_n(x)+(k_2-n\al\be)t(D)^2 P_n(x).
\end{equation*}
The desired difference or differential equation follows by substitution
of \eqref{17}( and some formal manipulations.

In general the paper is striking by its elegance. It is also
a purely mathematical paper, without any sign that the author
was educated as a theoretical physicist.
Furthermore, the paper has a certain terseness. The proofs of
many of the statements are silently left as short exercises for the reader.
\section{A predecessor, a missed class of polynomials
and a rediscovery}
\label{20}
Meixner \cite{Me3}
does not give any citations for his classes IV and V, so the reader
would think they are new. Still the class IV, i.e., the class of
orthogonal polynomials nowadays called Meixner polynomials, already
appears in 1931 in a paper by Truksa\footnote{See
\url{https://web.math.muni.cz/biografie/ladislav_truksa.html} (in Czech).
Ladislav Truksa, born 1891, studied actuarial science,
mathematics and physics
in Prague. His doctoral thesis in 1927 was on Legendre polynomials.
He habilitated in actuarial and mathematical statistics in 1931.
He played an important role in the early scientific and teaching
activities of the Department of Probability and Mathematical Statistics,
founded in 1952 within Charles University, Prague.} 
\cite{15}.
This author starts with the 3-parameter class of orthogonal polynomials
nowadays called Hahn polynomials \cite[Section 9.5]{1} which he can
trace back to Chebyshev \cite{16} (1873). In Part III Truksa \cite{15}
discusses all limit cases of the Hahn polynomials, as are now familiar from
the Askey scheme \cite[p.184]{1}. In particular, he obtains the
generalized Kummer polynomials \cite[p.192]{15}, which are the polynomials
nowadays called Meixner polynomials. He gives the weights, the hypergeometric
representation, the Rodrigues type formula, the three-term recurrence
relation, and the second order difference equation.
Curiously enough, both Truksa's paper \cite{15} and Meixner's paper
\cite{Me2} were reviewed in
{\em Jahrbuch \"uber die Fortschritte der Mathematik}
by W. Hahn (JFM 57.0413.01 and JFM 60.0293.01), but the reviewer did not
observe in his review of \cite{Me2} that the polynomials of class IV
already appeared in Truksa's paper. There are hardly any citations
of Truksa's paper in the literature. 
But Hahn mentions it in some of his JFM reviews.
He also mentions Truksa's and Meixner's paper in his important paper
\cite[p.32, footnote 1]{17}, but without giving details.

In addition to the five classes listed in the beginning of
Section \ref{21}, there should have been a sixth class as the result
of Meixner's classification:
\paragraph{Case VI}
Krawtchouk polynomials $K_n(x;p,N)$ ($n=0,1,\ldots,N$),
\cite[Section 9.11]{1}:
\begin{equation*}
\exp\left(x\log\left(\frac{1+(1-p^{-1})t}{1+t}\right)\right)
(1+t)^N=\sum_{n=0}^N \binom Nn K_n(x;p,N)\,t^n\quad(x=0,1,\ldots,N),
\end{equation*}
orthogonal on $\{0,1,\ldots,N\}$ with respect to the weights
$\binom Nx p^x (1-p)^{N-x}$.
\bLP
But Meixner cannot be blamed for this omission,
since he only looked for infinite systems $\{p_n\}_{n=0,1,2,\ldots}$
of orthogonal polynomials satisfying \eqref{2}. The inclusion of
Krawtchouk polynomials in Meixner's classification was probably first made
by Lancaster \cite{22} (1975).

The polynomials $P_n^{(\la)}(x;\phi)$ (case V of the classification)
were rediscovered by Pollaczek \cite{25} (1950). Quite remarkably,
he obtained them as limit cases of a three-parameter class of
orthogonal polynomials which he had introduced the year before and which
are nowadays known as Pollaczek polynomials. These latter polynomials
are of great importance but definitely outside the Askey scheme.
For quite a while the polynomials $P_n^{(\la)}(x;\phi)$ were only
attributed to Pollaczek in the literature. Around 1985, when the
Askey scheme was formalized, the name Meixner-Pollaczek polynomials
became common for these polynomials, see Andrews \& Askey \cite[p.41]{23}
for the motivation of this name.

As observed in
\cite[p.41]{23}, the Meixner, Krawtchouk and Meixner-Pollaczek polynomials
are essentially the same polynomials, but for different parameter ranges
and, in the Meixner-Pollaczek case, with a complex linear change of
independent variable. Indeed, from the expressions as hypergeometric
functions for the three classes, see \cite[Ch.~9]{1}, it follows that
\begin{align*}
K_n(x;p,N)&=M_n\left(x;-N,\frac p{p-1}\right)\quad(n=0,1,\ldots,N),\\
P_n^{(\la)}(x;\phi)&=\frac{(2\la)_n}{n!}\,e^{in\phi}
M_n(-\la-ix;2\la,e^{2i\phi}).
\end{align*}
We can group together the six cases of the Meixner classification
into the ``Meixner scheme'', a subscheme of the Askey scheme \cite[p.184]{1}.
The arrows indicate limit transitions, both for the polynomials and
for the orthogonality measures.
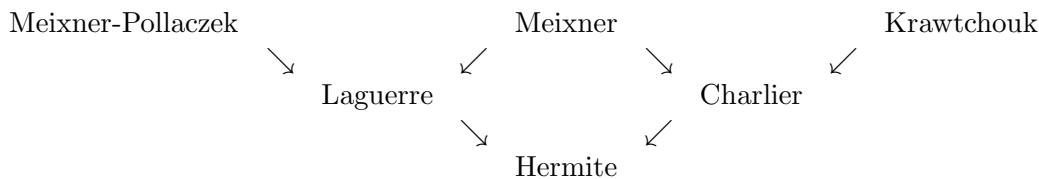
\[
\begin{matrix}
\mbox{Meixner-Pollaczek}&&&&{\rm Meixner}&&&&{\rm Krawtchouk}\\
&\searrow&&\swarrow&&\searrow&&\swarrow&\\
&&{\rm Laguerre}&&&&{\rm Charlier}&&\\
&&&\searrow&&\swarrow&&&\\
&&&&{\rm Hermite}&&&&
\end{matrix}
\]
One may also go directly from a family in the top row to Hermite in the bottom row by a limit.
\section{Follow-up of Meixner's 1934 paper}
\label{34}
We did some searches for texts which cite
Meixner's paper \cite{Me2} (1934), and in particular take up the new
classes of orthogonal polynomials introduced in \cite{Me2}.
Some early citations to \cite{Me2} concern the Charlier polynomials,
see for instance Szeg\"o's book \cite[\S2.81]{8} (1939) and Meixner's own
paper \cite{Me3}.
Early papers which continue work on the orthogonal polynomials
with generating function \eqref{2} are Geronimus \cite{9} (1938) and
Sheffer \cite{7} (1939). In particular, in \cite[Section~2]{7} a sequence
of polynomials $P_n(x)$ of degree $n$ satisfying
the equivalent properties \eqref{7} and \eqref{11} (slightly renormalized)
is called a set of type zero corresponding to the operator $t(D)$.
Sheffer gives there many further properties of such sequences (although
Meixner's operator $\La$ does not seem to figure in \cite{7}).

Clearly, there are many systems of monic polynomials $P_n(x)$ satisfying
\eqref{11} for given $t(D)$.
However, as observed by Sheffer \cite[Theorem 1.2]{7}, there
exists a unique such sequence which moreover satisfies $P_n(0)=0$ for
$n>0$. This is called the {\em basic sequence} associated with
$t(D)$\footnote{Clearly, a basic sequence cannot be a sequence of
orthogonal polynomials.}.
Steffensen \cite{10} (1941) independently introduced these basic sequences,
called {\em poweroids} by him.

Sheffer's sets of type zero \cite{7} got much attention in later papers,
without any reference to Meixner's paper \cite{Me2}. They became generally
known as {\em Sheffer polynomials}, see for instance~\cite{11}.
In particular, Sheffer polynomials became an important topic in the
so-called
{\em umbral calculus}, pioneered by Rota, Kahaner \& Odlyzko \cite{12},
see also the survey \cite{14} by Di Bucchianico.

The probability distributions corresponding to the five classes of
orthogonal polynomials resulting from Meixner's classification
\cite{Me2} got a lot of attention in papers on stochastics.
See surveys by Lai \cite{18} and Di Bucchianico \cite{14}.
These distributions also arise in the classification of the
so-called {\em natural exponential families} with quadratic variance
functions by Morris \cite{24}. See also \cite{g} on stochastic
aspects of Meixner polynomials.

Meixner and Meixner-Pollaczek polynomials have group theoretic
interpretations as matrix elements of discrete series representations
of SL(2,$\RR$) (or of its universal covering), see
Basu \& Wolf \cite[Section 3]{26} and
Koornwinder \cite[Section 7]{27}.
Foata \& Labelle \cite{35} gave combinatorial models for the
Meixner polynomials.

It often happens with mathematical objects named after a person that
a generalized object also keeps this name, although the generalization
has been made by someone else. Thus happened also with Meixner polynomials.
We meet for instance:
\begin{itemize}
\item
{\em $q$-Meixner polynomials}, \cite[Section 14.13]{1}.
These $q$-analogues of Meixner polynomials first occurred in Hahn's
classification \cite{17}, but not yet bearing this name.
\item
{\em Multivariable Meixner polynomials}. These were introduced by
Tratnik \cite{28}. They involve (non-straightforward) products of
one-variable Meixner polynomials. They could already have been,
but were not, obtained as limits of the multivariable Hahn polynomials
considered by Karlin \& McGregor \cite{29}. See Gasper \& Rahman
\cite{37} for the multivariable $q$-case.
In principle, another kind of multivariable Meixner polynomials,
associated with root system $BC_n$, can be obtained as a limit case
of the $BC_n$ Hahn polynomials considered by van Diejen \& Stokman
\cite[Section 5]{30}. Possibly the papers by Chikuse \cite{32}
and Bryc \& Letac \cite{31}, written from a probabilistic point of view,
and the paper by Shibukawa \cite{34} are related.
\end{itemize}
One also meets multiple Meixner polynomials, exceptional Meixner
polynomials, $d$-orthogonal Meixner polynomials,
matrix-valued Meixner polynomials and Meixner-Sobolev
orthogonal polynomials in the literature.
In stochastics one finds Meixner distributions, Meixner ensembles
and Meixner processes. On September 1, 2016 MathSciNet listed
188 publications having the word {\em Meixner} in its title.
Similarly zbMATH found 182 documents.

As a concluding remark, the irony of history is very visible
in connection with Meixner's 1934 paper \cite{Me2}. Indeed,
Meixner introduced the Sheffer polynomials,
but they are called after Sheffer.
Meixner formulated the Favard theorem, but it is called after Favard.
Meixner gave the Meixner polynomials, but they should be called
after Truksa.
Meixner gave the Meixner-Pollaczek polynomials, but they were almost
called just Pollaczek polynomials.
\subsection*{Acknowledgements}
We are very grateful to Dr.~Michael Meixner (Aachen),
who provided one of us (PB)
with many documents about his father Josef Meixner.
We also express our thanks to
Prof.~Janine Splettstoesser (G\"oteborg), who passed to us some
recollections by people who were around in Aachen as assistant or
student in the days of Meixner's professorship.
Thanks also to Prof.~Ernst G\"orlich and Prof.~Clemens Markett
(both from Aachen),
who helped us with literature (the latter also kindly pointed out
some errata in a previous version).
We thank Prof.~Alessandro Di Bucchianico (Eindhoven)
for helpful comments about Sheffer polynomials.
Finally we thank \mbox{Andrius} Kulikauskas, who observed in the previous version
an unforgivable error in our expression of the weights for which the Meixner
polynomials are orthogonal.
\renewcommand{\refname}{Selected Publications by Meixner\footnote{\cite{Me1}--\cite{Me16} are papers on orthogonal
polynomials and hypergeometric functions,
\cite{Me7}--\cite{Me15b} are papers on
spheroidal wave functions and Mathieu functions, and
\cite{B}--\cite{MeSchWo} are books and Handbuch articles.}
}
\let\oldthebibliography=\thebibliography
\let\oldendthebibliography=\endthebibliography
\renewenvironment{thebibliography}[1]{%
    \oldthebibliography{#1}%
    \setcounter{enumiv}{32}%
}{\oldendthebibliography}
\renewcommand{\refname}{Other references}

\vskip0.5cm
\begin{footnotesize}
\noindent
P. L. Butzer, Lehrstuhl A f\"ur Mathematik, RWTH Aachen,
52056 Aachen, Germany;
\sPP
email: {\tt butzer@rwth-aachen.de}
\bLP
T. H. Koornwinder, Korteweg-de Vries Institute, University of
Amsterdam,\sPP
P.O.\ Box 94248, 1090 GE Amsterdam, The Netherlands;
\sPP
email: \texttt{thkmath@xs4all.nl}
\end{footnotesize}

\end{document}